\documentclass[11pt, reqno]{amsart}
\usepackage{latexsym}
\usepackage{amsmath, amssymb, amsfonts}
\author{Chengyan Liu}
\address{7003 Brocton Ct., Virginia, VA 22150}
\email{geismn@bigfoot.com}
\title{Sieve Method and Landau Problems}
\begin{document}
\theoremstyle{plain}
\newtheorem{goman}{Theorem}
\newtheorem{nason}{Lemma}
\newtheorem{myson}[nason]{Lemma}
\newtheorem{meree}[nason]{Lemma}
\newtheorem{mersn}[nason]{Lemma}
\newtheorem{cnver}[nason]{Lemma}
\newtheorem{panpn}[nason]{Lemma}
\newtheorem{dilea}[nason]{Lemma}
\theoremstyle{remark}
\newtheorem*{erfem}{Feature M}
\newtheorem*{gcsiu}{Landau Problem 1}
\newtheorem*{tpsiu}{Landau Problem 2}
\newtheorem*{sqsiu}{Landau Problem 3}
\newtheorem*{spqsu}{Landau Problem 4}
\newtheorem*{fnalr}{Remark}
\begin{abstract}
We solve Landau's four {\em unattackable problems\/}, including Goldbach 
Conjecture and Twin Prime Conjecture through sieve method.
\end{abstract}
\maketitle
At the 1912 Fifth International Congress of Mathematicians (Cambridge), 
Landau mentioned four {\em unattackable problems\/}: Goldbach Conjecture; 
Twin Prime Conjecture; for each $n$, there is a prime $p$ such that 
$n^{2} < p < (n + 1)^{2}$; there are infinitely many primes $p$ of the form 
$p = n^{2} + 1$. As Erd\H os~\cite[preface]{gur:pnt} maintained five years ago, 
we still do not have a satisfactory solution for each of them. 

For large $N$, let $\mathcal{A} = \{ m |~m \in \mathbb{Z}, m \leq N \}$, 
$\mathcal{P} = \{ p |~p \in \mathcal{A}, p~\text{prime} \}$, and
\begin{equation}\label{ari:fnn}m = g(p)~,\end{equation}such that 
$m \in \mathcal{A}$, $m \not \in \mathcal{P}_{0}$, 
$p \in \mathcal{P}_{z} = \{ p |~p \in \mathcal{P}, 2 < p \leq z < N \}$, 
$\mathcal{P}_{0}$ a desired set. A rudimentary sieve associated with 
Eratosthenes employs\begin{quotation}\begin{erfem}Screening out 
$m \not \in \mathcal{P}$ through~\eqref{ari:fnn} by 
$p \in \mathcal{P}_{\surd{N}}$.\end{erfem}\end{quotation}However, 
since its residue term is proportionate to 
card$\{ q : p~|~q, p \in \mathcal{P}_{\surd{N}}\}$, we could not use it for 
practical purpose. 

Yet this is not our concern here since there are many amendments to 
improve its estimation for both main and residue terms. Our concern here is 
{\em Feature M\/} which virtually all of the revision retain. It leads to a lengthy 
and often difficult elimination process according to how many prime factors 
$m$ have. In other words, it limits our focus on excluding composite numbers. 
For Goldbach Conjecture, Chen~\cite{chj:ept:ine, chj:ept:ngl} was stalled for 
$m$ has at most two prime factors. Nevertheless, we could unlock our 
cogitation by adding an inclusion scheme. We rewrite~\eqref{ari:fnn} as
\begin{equation}\label{ari:fss}p = g~(\!\!\prod_{q~\in~\mathcal{P}_{\surd{N}}}\!\! q~)~,
\end{equation}where $p \in \mathcal{P}$. With this function, we collect certain 
primes $p$ first, usually through an arithmetic progression, and then focus 
on discharging them if $p \not \in \mathcal{P}_{0}$. It enhances sieve methods 
and their revisions with flexibility.

The latest development of sieve methods is a Bombieri-Vinogradov
~\cite{boe:lar, via:den:lun, via:rec:den:lun} type theorem. It treats residue 
term in the following fashion,\begin{align}
& \sum_{d \leq D}\underset{y \leq x}{max}\underset{(l, d) = 1}{max}
    \biggl|\sum_{\substack{n \leq y \\ n~\equiv~l \!\!\!\!\!\pmod{d}}}
                \Lambda(n) - \frac{y}{\varphi(d)} \biggr| \notag \\
   = & \sum_{d \leq D}\underset{y \leq x}{max}\underset{(l, d) = 1}{max}
          \biggl|\frac{1}{\varphi(d)}\sum_{\chi}\Bar{\chi}(l)\psi^{\prime}(y, \chi)\biggr| \notag \\
   \leq & \sum_{d \leq D}\underset{y \leq x}{max}\underset{(l, d) = 1}{max}
               \frac{1}{\varphi(d)}\sum_{\chi}|\psi^{\prime}(s, \chi)| \label{bom:vin}\\
   \ll & \frac{x}{(\log x)^{A}}~, \notag\end{align}
for $D = x^{\frac{1}{2}}(\log x)^{- B}$, where $B = A + \eta$, $\eta \geq 2$ a constant, 
$A$ an arbitrary given number, 
$\psi^{\prime}(y, \chi) = \psi(y, \chi) = \sum_{n \leq y} \chi(n) \Lambda(n)$,
except when $\chi = \chi_{0}$, $\psi^{\prime}(y, \chi) = \psi(y, \chi) - y$, $\chi$ a 
Dirichlet character modulus $d$, $\Lambda(n)$ the Manglodt function, 
$\varphi(d)$ the Euler function.

Estimation step occurred at~\eqref{bom:vin} is our other concern. For
an {exclusion-inclusion} process, the step means that we only take an 
{\em exclusion\/} process without sufficient compensation. At best, on 
Riemann Hypothesis, we could have $\eta = 2$ but leave a gap 
$D = x^{\frac{1}{2}}(\log x)^{- A - 2}$. To close the gap, we deal with
\begin{equation}\label{caa:eai}
\sum_{d \leq D} \frac{\mu(d)}{\varphi(d)} \sum_{\chi \not = \chi_{0}} \Bar{\chi}(a) \psi(x, \chi)~.
\end{equation}In~\eqref{caa:eai}, we take explicit formula
\footnote{See~\cite[\S 19]{dah:mnt}.} for $\psi(x, \chi)$
\begin{align}
    \psi(x, \chi) = &~- \sum_{|Im(\rho)| < T} \frac{x^{\rho}}{\rho} + \sum_{|Im(\rho)| < 1} \frac{1}{\rho}
                                + c_{1} \frac{x (\log dx)^{2}}{T} \notag \\
=  &~\psi_{\rho}(d, x, T)~,\label{caa:etu}\end{align}where $\rho$ are nontrivial zeros 
of Dirichlet $L$-function $L(s, \chi)$. Since there is no term in~\eqref{caa:etu}  
directly related to $\chi$, taking out $\psi_{\rho}(d, x, T)$ in~\eqref{caa:eai} we have\begin{align}
& \sum_{d \leq D} \frac{\mu(d)}{\varphi(d)}~\psi_{\rho}(d, x, T) 
                                \sum_{\chi \not = \chi_{0}} \Bar{\chi}(a) \notag \\
                        =~ & \sum_{d \leq D} \frac{\mu(d)}{\varphi(d)}~\psi_{\rho}(d, x, T) 
                                \biggl(\varphi_{a, 1}(d) - 1 \biggr)~,\label{caa:etr}\end{align}where 
$\varphi_{a, 1}(d) = \varphi(d)$ if $a \equiv 1 \pmod{d}$ and $\varphi_{a, 1}(d) = 0$ otherwise. 
Now, we could take absolute value for~\eqref{caa:etr} by using following lemmata
\begin{nason}\label{snn:eai}\begin{equation*}\sum_{n \leq x} \frac{1}{\varphi(n)} \ll \log x~.
\end{equation*}\end{nason}
\begin{proof}This is Theorem A.17 of Nathanson~\cite[p. 316]{nam:ant:tcb}.\end{proof}
\begin{myson}\label{snn:etu}All nontrivial zeros of Dirichlet $L$-function $L(s, \chi)$ lie on 
line $Re(s) = 1/2$.\end{myson}\begin{proof}This is Theorem 2 of Liu~\cite{lic:rhy}.
\end{proof}

By {\bf Lemma\/}~\ref{snn:eai} and {\bf Lemma\/}~\ref{snn:etu}, we can take
$T = x^{\frac{1}{2}}$ and $Re(\rho) = 1/2$ in~\eqref{caa:etu} to get
\begin{equation}\label{caa:efo}
\biggl|\sum_{d \leq D} \frac{\mu(d)}{\varphi(d)} 
                       \sum_{\chi \not = \chi_{0}} \Bar{\chi}(a) \psi(x, \chi)\biggr| 
              \leq O(x^{\frac{1}{2}} (\log x)^{3})~,\end{equation}
where $D = x^{\frac{1}{2}}$. To implement a sieve with~\eqref{ari:fss} 
and~\eqref{caa:efo}, we need 
\begin{meree}[Mertens]\label{muu:eai}For $x \geq 1$,\begin{equation*}
\sum_{p \leq x} \frac{\log p}{p} = \log x + O(1)~.\end{equation*}\end{meree}
\begin{proof}This is Theorem 6.6 of Nathanson~\cite[p. 160]{nam:ant:tcb}.\end{proof}
\begin{mersn}[Mertens]\label{muu:etu}For $x \geq 2$,\begin{equation*}
\prod_{p \leq x} \biggl(1 - \frac{1}{p}\biggr)^{- 1} = e^{\gamma} \log x + O(1)~,\end{equation*}
where $\gamma$ is Euler constant.\end{mersn}
\begin{proof}This is Theorem 6.8 of Nathanson~\cite[p. 165]{nam:ant:tcb}.\end{proof}
\begin{cnver}\label{muu:etr}\begin{gather*}
\prod_{q} \biggl(1 - \frac{1}{\varphi(q^{\theta})}\biggr)\biggl(1 - \frac{1}{q}\biggr)^{- 1}~, \\
\intertext{converge, and}
\prod_{q \geq w} \biggl(1- \frac{1}{\varphi(q^{\theta})}\biggr)\biggl(1 - \frac{1}{q}\biggr)^{- 1} 
                                 = 1 + O\biggl(\frac{1}{\log w}\biggr)~,\end{gather*}where $\theta$ is
a given number.\end{cnver}
\begin{proof}Applying {\bf Lemma\/}~\ref{muu:eai}, this is Corollary 2 for Lemma 7 of 
Pan and Pan~\cite[\S 7, p. 163]{pdn:pbn:tgc}.\end{proof}
\begin{panpn}\label{muu:pun}For $2 \leq w \leq z$ and a given $\theta$, we have
\begin{equation*}\prod_{w \leq q < z} \biggl(1 - \frac{1}{\varphi(q^{\theta})}\biggr)
            = \frac{\log z}{\log w}\Biggl(1 + O\biggl(\frac{1}{\log w}\biggr)\Biggr)~.\end{equation*}
\end{panpn}\begin{proof}Applying {\bf Lemma\/}~\ref{muu:eai}, this is Corollary 3 for
Theorem 2 of Pan and Pan~\cite[\S 7, p. 164]{pdn:pbn:tgc}.\end{proof}
\begin{dilea}[Dirichlet]\label{lai:dnt}
\begin{equation*}\pi (l, d; x) = \! \sum_{\substack{p~\leq~x \\ p~\equiv~l \!\!\!\!\! \pmod d}} \!\!\! 1 = 
                                  \frac{x}{\varphi(d) \log x} + O(x^{1 \over 2} \log x)~,\end{equation*}
where $p$ is a prime, $d \leq x$, $(l, d) = 1$, $\varphi(d)$ the Euler function.\end{dilea}
\begin{proof}By Theorem 3.5.1 of Gelfond and Linnik~\cite[\S 3.5, p. 72]{gea:ent:liy}, we have
\begin{equation*}
\pi (l, d; x) = \frac{x}{\varphi(d) \log x} + o\biggl(\frac{x}{\log x}\biggr)~.\end{equation*}
Applying {\bf Lemma\/}~\ref{snn:etu} and the remark of Davenport~\cite[\S 20]{dah:mnt}, 
we have our result.\end{proof}

\begin{gcsiu} This is Goldbach Conjecture. We take $p = N - m$ as our function,
where $N \equiv 0 \pmod{2}$, $q \in \mathcal{P}_{\surd{N}}$, $q~|~m$, $q \nmid N$,
$p \in \mathcal{P}$, and try to find out primes $p$ belongs to arithmetic progression
$p \equiv N \pmod{q}$. To get an answer, we need to remove terms 
$r_{q} - q\nu$, where $r_{q} \equiv N \pmod{q}$, $1 \leq r_{q} < q$, $\nu > 1$, from 
the progression. Applying {\bf Lemma\/}~\ref{muu:etu}, {\bf Lemma\/}~\ref{muu:etr}, 
{\bf Lemma\/}~\ref{muu:pun}, and {\bf Lemma\/}~\ref{lai:dnt} with $\theta = 1$, we
have (main term)\begin{align}
    & \prod_{\substack{q \nmid N \\ q \leq \surd{N}}} 
        \biggl(1 - \frac{1}{\varphi(q)}\biggr) \frac{N}{\log N}\notag \\
=  & \prod_{\substack{q \nmid N \\ q < N}} \biggl(1 - \frac{1}{\varphi(q)}\biggr)
        \prod_{\substack{q \nmid N \\ \surd{N} < q < N}} \frac{q - 1}{q - 2} \frac{N}{\log N}\notag \\
= &  \prod_{q} \biggl(1 - \frac{1}{\varphi^{2}(q)}\biggr)
        \prod_{\substack{q | N \\ q < N}} \frac{q - 1}{q - 2}
        \prod_{\substack{q \nmid N \\ \surd{N} < q < N}} \frac{q - 1}{q - 2} 
        \frac{e^{- \gamma}N}{(\log N)^{2}} \Biggl(1 + O\biggl(\frac{1}{\log N}\biggr)\!\!\Biggr)\notag \\ 
= &  \prod_{q} \biggl(1 - \frac{1}{\varphi^{2}(q)}\biggr)
        \prod_{\substack{q | N \\ q < N}} \frac{q - 1}{q - 2}
        \frac{2 e^{- \gamma} N}{(\log N)^{2}} \Biggl(1 + O\bigg(\frac{1}{\log N}\biggr)\!\!\Biggr)~.
\label{gaa:etn}\end{align}The above equation might also eliminate $p = N - q$ for 
$q \leq \surd{N}$. At most, it costs us $O(\surd{N})$. By~\eqref{caa:efo}, we have 
residue term\begin{equation}\label{gaa:ele}O(N^{\frac{1}{2}} (\log N)^{2})~.\end{equation}\end{gcsiu}
\begin{tpsiu}This is Twin Prime Conjecture. We take $p - 2 = m$ as our function,
where $p \in \mathcal{P}$, $q | m$, $q \in \mathcal{P}_{\surd{N}}$, and try to find 
out primes $p$ in arithmetic progression $p \equiv 2 \pmod{q}$. For this purpose,
we can proceed just as the above. To find out number of prime pairs such that 
$p_{1} - 2^{\alpha}K = p_{2}$, where $K \equiv 1 \pmod{2}$, $\alpha \geq 1$, 
$p_{j} \in \mathcal{P}$, $j = 1, 2$, we take $p - 2^{\alpha}K = m$ as our function, 
where $p \in \mathcal{P}$, $q | m$, $q \in \mathcal{P}_{\surd{N}}$, $q \nmid K$. 
As before, we want to eliminate terms $r + q\nu$, where $r \equiv 2^{\alpha} K \pmod{q}$, 
$1 \leq r < q$, $\nu > 1$, from arithmetic progression $p \equiv 2^{\alpha}K \pmod{q}$. 
If $p = 2^{\alpha} K + q$, $q \leq \surd{N}$, we might loss at most 
$O(\surd{N})$ primes. Otherwise, we would end up with exactly~\eqref{gaa:etn} 
and~\eqref{gaa:ele} for main and residue term respectively with $q | K$ substitute for 
$q | N$ in~\eqref{gaa:etn}.\end{tpsiu}
\begin{sqsiu}We treat $N \equiv 0 \pmod{2}$ first. Assume $K \equiv 1 \pmod{2}$, 
$N = 2^{\alpha} K$, we use equation $p = N^{2} + r$ to find out prime $q < N$ such that 
$q^{2} | p - r$, $q | K$, where $1 \leq r \leq 2N$, 
$p \in \mathcal{P}_{1} = \{ p |~N^{2} < p < (N + 1)^{2}, p~\text{prime}\}$, $r$, $N^{2}$ 
and $p$ pairwise coprime. It is obvious if $p \in \mathcal{P}_{1}$ then it belongs to
arithmetic progression $p \equiv r \pmod{N^{2}}$, $(r, N^{2}) = 1$. Hence, 
number of primes (main term) such that $p \in \mathcal{P}_{1}$ equals to
\begin{equation}\label{gaa:sin}\frac{(N + 1)^{2}}{\varphi(N^{2}) \log (N + 1)^{2}}~,\end{equation}
by {\bf Lemma\/}~\ref{lai:dnt}. Yet, there are\begin{equation}\label{gaa:cos}
\frac{1}{2} \prod_{\substack{q | K \\ q < N}} \biggl(1 - \frac{1}{q}\biggr) 2N~,\end{equation}
$r$ satisfy~\eqref{gaa:sin}. We put~\eqref{gaa:sin} and~\eqref{gaa:cos} together
\begin{align*}
    &  \frac{1}{\varphi(N^{2})} \prod_{\substack{q | K \\ q < N}} \biggl(1 - \frac{1}{q}\biggr) 
              \frac{N (N + 1)^{2}}{\log (N + 1)^{2}} \\
= & \frac{N}{\log N} \Biggl(1 + O\biggl(\frac{1}{N}\biggr)\Biggr)~.\end{align*} 
For residue term, we have\begin{align*}
    &~ \prod_{\substack{q | K \\ q < N}} \biggl(1 - \frac{1}{q}\biggr) N \frac{2N}{\varphi(N^{2})} \\
=  &~ \frac{1}{2^{2 \alpha - 1}} \prod_{\substack{q | K \\ q < N}} \frac{1}{q^{2}} ~2 N^{2} \\
=  &~ 4~.\end{align*}We can treat $N + 1 \equiv 0\!\pmod{2}$ the same way by substituting 
$p = N^{2} + r$ with $p = (N + 1)^{2} - r$.\end{sqsiu}
\begin{spqsu}For large $N$, we use function $n^{2} + 1 = m$ to find out primes 
$p$ of the form $n^{2} + 1$, where $n < \surd{N}$. For $m$ to be a candidate, it 
suffices $n \equiv 0 \pmod{2}$ and $m \equiv 1 \pmod{4}$. If prime $q$ satisfies 
$q | m$, it is necessary that $q \equiv 1 \pmod{4}$. For a fixed $q \equiv 1 \pmod{4}$, 
number of $n^{2}$ (main term) belong to arithmetic progression 
$n^{2} \equiv - 1 \pmod{q}$ equals to\begin{equation*}
\frac{1}{\varphi(q)} \prod_{\substack{q~\equiv~1\!\!\!\pmod{4} \\ q \leq \surd{N}}} 
                                     \biggl(1 - \frac{1}{q}\biggr) \surd{N}~,\end{equation*}
because $q \nmid n$. To remove terms $n^{2} = q \nu - 1$, $\nu > 1$ from the 
progression, we have
\begin{align*}
    & \!\!\!\prod_{\substack{q~\equiv~1\!\!\!\pmod{4} \\ q \leq \surd{N}}}\!\!\! 
        \biggl(1 - \frac{1}{\varphi(q)}\biggr) 
        \!\!\!\prod_{\substack{q~\equiv~1\!\!\!\pmod{4} \\ q \leq \surd{N}}}\!\!\!
        \biggl(1 - \frac{1}{q}\biggr) \sqrt{N} \\
=  & \!\!\!\prod_{q~\equiv~3\!\!\!\!\pmod{4}}\!\!\!
        \biggl(1 + \frac{1}{\varphi(q)}\biggr)
        \!\!\!\prod_{q~\equiv~1\!\!\!\!\pmod{4}}\!\!\! 
        \biggl(1 - \frac{1}{\varphi(q)}\biggr) \frac{2 e^{- \gamma} \surd{N}}{\log N}
        \Biggl(1 + O\biggl(\frac{1}{\log N}\biggr)\!\!\Biggr)~.\end{align*}
Notice that it might also eliminate $q - 1 = n^{2}$ for $n < \sqrt[4]{N}$, $q \leq \surd{N}$. 
With this, we might loss at most $O(\sqrt[4]{N})$ primes. For residue term, 
applying~\eqref{caa:efo} we get \begin{equation*}O(\sqrt[4]{N}(\log N)^{2})~.\end{equation*}
since total terms for $n^{2} < N$ is $\surd{N}$.\end{spqsu} 
%

%
\end{document}